\documentclass[11pt]{article}
\usepackage[all]{xy}

\usepackage{latexsym}
\usepackage{amsfonts}
\usepackage{amssymb}
\usepackage{mathrsfs}
\newcommand{\be}{\begin{equation}}
\newcommand{\ee}{\end{equation}}
\newcommand{\bea}{\begin{eqnarray}}
\newcommand{\eea}{\end{eqnarray}}
\newcommand{\bean}{\begin{eqnarray*}}
\newcommand{\eean}{\end{eqnarray*}}
\newcommand{\brray}{\begin{array}}
\newcommand{\erray}{\end{array}}

\newcommand{\newsection}[1]{\setcounter{equation}{0}
\setcounter{dfn}{0}
\section{#1}}

\newtheorem{dfn}{Definition}[section]
\newtheorem{thm}[dfn]{Theorem}
\newtheorem{lmma}[dfn]{Lemma}
\newtheorem{ppsn}[dfn]{Proposition}
\newtheorem{crlre}[dfn]{Corollary}
\newtheorem{xmpl}[dfn]{Example}
\newtheorem{rmrk}[dfn]{Remark}

\newcommand{\bdfn}{\begin{dfn}}
\newcommand{\bthm}{\begin{thm}}
\newcommand{\blmma}{\begin{lmma}}
\newcommand{\bppsn}{\begin{ppsn}}
\newcommand{\bcrlre}{\begin{crlre}}
\newcommand{\bxmpl}{\begin{xmpl}}
\newcommand{\brmrk}{\begin{rmrk}}

\newcommand{\edfn}{\end{dfn}}
\newcommand{\ethm}{\end{thm}}
\newcommand{\elmma}{\end{lmma}}
\newcommand{\eppsn}{\end{ppsn}}
\newcommand{\ecrlre}{\end{crlre}}
\newcommand{\exmpl}{\end{xmpl}}
\newcommand{\ermrk}{\end{rmrk}}

\newcommand{\IC}{\mathbb{C}}
\newcommand{\bbc}{\mathbb{C}}
\newcommand{\IZ}{\mathbb{Z}}
\newcommand{\bbz}{\mathbb{Z}}
\newcommand{\IN}{\mathbb{N}}\newcommand {\bbn}{\mathbb{N}}

\newcommand{\cla}{\mathcal{A}}
\newcommand{\clb}{\mathcal{B}}
\newcommand{\clh}{\mathcal{H}}
\newcommand{\cli}{\mathcal{I}}
\newcommand{\clk}{\mathcal{K}}
\newcommand{\cll}{\mathcal{L}}
\newcommand{\clq}{\mathcal{Q}}

\newcommand{\prf}{\noindent{\it Proof\/}: }
\newcommand{\ots}{\otimes}
\newcommand{\raro}{\rightarrow}

\newcommand{\seq}{\subseteq}

\newcommand{\nn}{\nonumber}

\def \qed { \mbox{}\hfill
$\Box$\vspace{1ex}}

\begin{document}

\author{\textsc{Partha Sarathi Chakraborty} 
 and 
\textsc{Arupkumar Pal}}
\title{Spectral triples and  
associated Connes-de Rham complex for the quantum $SU(2)$ and the quantum sphere}

\maketitle


\begin{abstract}
  In this article, we construct  spectral triples for the
$C^*$-algebra of continuous functions on the quantum $SU(2)$ group
and the quantum sphere.
There has been various approaches towards building a calculus on quantum
spaces, but there seems to be very few instances of computations outlined in chapter~6,
\cite{C7}. We give detailed computations of the associated Connes-de Rham 
complex and the space of $L_2$-forms.
\end{abstract}
{\bf AMS Subject Classification No.:} {\large 58}B{\large 34},
{\large 81}R{\large 50}, {\large 46}L{\large 87}\\
{\bf Keywords.} Spectral triples, exterior complex.


\section{Introduction}

Given a noncommutative space, there is no general method for 
constructing a spectral triple on it. 
Even though there are general results asserting the existence of 
enough unbounded Kasparov modules (\cite{BJ}),
in concrete examples, it is often difficult to carry out this prescription.
In \cite{CP}, the authors characterized all spectral triples for 
the $C^*$-algebra $\cla$ of continuous functions on $SU_q(2)$
represented on its $L_2$-space, assuming equivariance under 
the (co-)action of the group itself.
In the present article, we take the more standard representation
of $\cla$ on $\clh=L_2(\bbn)\otimes L_2(\bbz)$ (see (\ref{pi}) below),
and impose equivariance condition under the action of the group $S^1\times S^1$.
Employing similar technques as in \cite{CP}, we arrive at a spectral triple 
of dimension 2. One advantage of this triple is that it is relatively 
easy to compute the associated Connes-de Rham complex, which we 
give in section~3. 
This complex is supported on $\{0,1\}$, and thus captures
the topological dimension, which can be seen to be 1 from the 
following well-known exact sequence
\be
0\longrightarrow \clk\otimes C(S^1) \longrightarrow  \cla
 \longrightarrow C(S^1) \longrightarrow 0.
\ee
The complex of square integrable forms were introduced by
Frolich et.\ al in \cite{FGR}. We also present calculations of these 
$L_2$-forms for this spectral triple.

In the last section, we briefly indicate how to carry out a similar construction 
of a spectral triple and the associated calculus for the quantum sphere $S^2_{qc}$.

Let us start with a brief description of the $C^*$-algebra 
of continuous functions on the quantum $SU(2)$, 
to be denoted by $\cla$.
This is the  canonical $C^*$-algebra generated by two
elements $\alpha$ and $\beta$ satisfying the following relations:
\be
\alpha^*\alpha+\beta^*\beta=I,\;\alpha\alpha^*
+q^2\beta\beta^*=I,\;
\alpha\beta-q\beta\alpha=0,\;\alpha\beta^*-q\beta^*\alpha=0,\;
\beta^*\beta=\beta\beta^*.
\ee

The $C^*$-algebra $\cla$ can be described more concretely as
follows. Let $\{e_i\}_{i\geq0}$ and $\{e_i\}_{i\in\IZ}$ be the
canonical orthonormal bases for $L_2(\IN)$ and $L_2(\IZ)$ respectively.
We denote by the same symbol $N$ the operator $e_k\mapsto ke_k$,
$k\geq0$, on $L_2(\IN)$ and $e_k\mapsto ke_k$, $k\in\IZ$, on $L_2(\IZ)$.
Similarly, denote by the same symbol $\ell$ the operator $e_k\mapsto
e_{k-1}$, $k\geq1$, $e_0\mapsto0$ on $L_2(\IN)$ and the operator
$e_k\mapsto e_{k-1}$, $k\in\IZ$ on $L_2(\IZ)$. Now take $\clh$ to be
the Hilbert space $L_2(\IN)\ots L_2(\IZ)$, and define $\pi$ to be the
following representation of $\cla$ on $\clh$:
\be\label{pi}
\pi(\alpha)=\ell\sqrt{I-q^{2N}}\ots I,\;\;\;\pi(\beta)=q^N\ots \ell.
\ee
Then $\pi$ is a faithful representation of $\cla$, so that one
can identify $\cla$ with the $C^*$-subalgebra of $\cll(\clh)$
generated by $\pi(\alpha)$ and $\pi(\beta)$. Image of $\pi$ contains 
$\clk \otimes C(S^1)$ as an ideal with $C(S^1)$ as the quotient algebra, 
that is we have a useful short exact sequence
\be
\label{100}
0\longrightarrow \clk\otimes C(S^1) \stackrel i\longrightarrow  \cla
\stackrel\sigma \longrightarrow C(S^1) \longrightarrow 0.
\ee
We will denote by $\cla_f$ the *-subalgebra of $\cla$ 
generated by $\alpha$ and $\beta$. Let 
\[
\alpha_i\beta^j\beta^{*k}:=\cases{\alpha^i\beta^j\beta^{*k} & if $i\geq 0$,\cr
                                   (\alpha^*)^{-i}\beta^j\beta^{*k} & if $i<0$.}
\]
Then $\{\alpha_i\beta^j\beta^{*k} : i\in\bbz, j,k\in\bbn\}$ is a basis 
for $\cla_f$.
The Haar state $h$ on $\cla$ is given by, 
\[
h :a \mapsto (1-q^2)\sum_{i=0}^\infty q^{2i} \langle e_{i0},a e_{i0} \rangle. 
\]

\brmrk \label{interpret}\rm 
The representation $\pi$ admits a nice interpretation. Let $M$
  be a compact topological manifold and $E$, a Hermitian vector bundle
  on $M$. Let $\Gamma(M,E)$ be the space of continuous sections. Then
  $\Gamma(M,E)$ is a finitely generated projective $C(M)$ module.
  Define an inner product on $\Gamma(M,E)$ as
  \[
  \langle s_1,s_2 \rangle := \int {( s_1(m),s_2 (m))}_m d \nu (m),
  \]
  where $\nu$ is a smooth measure on $M$ and ${(\cdot,\cdot)}_m$
  is the inner product on the fibre on $m$. Let $\clh_E$ be the
  Hilbert space completion of $\Gamma(M,E)$. Then we have a natural
  representation of $C(M) $ in $ \cll(\clh_E)$. The same program can
  be carried out in the noncommutative context also. Let $\clb$ be a
  $C^*$-algebra and $E$ a Hilbert $\clb$-module with its $\clb$
   valued inner product ${\langle \cdot,\cdot \rangle}_\clb $. Let
  $\tau$ be a state on $\clb$. Consider the inner product on $E$ given
  by 
$\langle e_1, e_2 \rangle = \tau ( {\langle e_1, e_2 \rangle}_\clb )$. 
If we denote by $\clh_E$ the Hilbert space
  completion of $E$, then we get a natural representation of $\clb$ in
  $\cll(\clh_E)$. Now in the context of $SU_q(2)$, let 
$p= | e_0\rangle \langle e_0| \otimes I \in \cla$. 
Then it is easy to verify that $\clh_E=l^2(\bbn) \otimes l^2(\IZ)$ 
for $E=\cla p$  with its natural left Hilbert $\cla$-module structure.
  Moreover, the associated representation is nothing but the 
  representation of $\cla$ described above. Also, viewed this way, 
one can think of the representation of $\cla$ on $L_2(h)$ given in 
$\cite{CP}$ as being a countable direct sum of representations 
each of which look like $\pi$ (just think of $\cla$ as $\oplus \cla p_i$
where $p_i=| e_i\rangle \langle e_i| \otimes I$).
\ermrk

\newsection{$S^1\times S^1$-equivariant spectral triples}
The group $G=S^1\times S^1$ has the following action on $\cla$:
\be
\tau_{z,w}:\cases{\alpha\mapsto z\alpha &\cr
                  \beta\mapsto w\beta &}
\ee
Let $U$ be the following representation of $G$ on $\clh$:
$U_{z,w}=z^N\otimes w^N$. Then for any $a\in\cla$,
one has 
$\pi(\tau_{z,w}(a))= U_{z,w}^*\pi(a)U_{z,w}$, i.e.\ the action 
$\tau$ is implemented through this representation $U$ of $G$.
A self-adjoint operator with discrete spectrum equivariant under
this G-action must be of the form
\begin{equation}
D: e_{ij}\mapsto d_{ij}e_{ij}.
\end{equation}
It is easy to see that if $D$ is such an operator, then 
$[D,\alpha]$ and $[D,\beta]$ are given by
\begin{eqnarray} \label{eq:comm1}
[D,\alpha]e_{ij}  &=& (d_{i-1,j}-d_{ij})\sqrt{1-q^{2i}}\,e_{i-1,j},\\
{}[D,\beta]e_{ij} &=& (d_{i,j-1}-d_{ij})q^{i}\,e_{i,j-1}.\label{eq:comm1a}
\end{eqnarray}


Employing arguments very similar to those used in the proofs of 
propositions~3.1 and 3.2 in \cite{CP}, we now get the following results.

\bppsn \label{condn1}
Let $D$ be an operator of the form  $e_{ij}\mapsto d_{ij}e_{ij}$.
Then $[D,a]$ is bounded for all $a\in\cla_f$ if and only if $d_{ij}$'s satisfy
the following two conditions:
\begin{eqnarray}\label{bnd1}
|d_{i-1,j}-d_{ij}| &=&O(1),\\
 \label{bnd3}
|d_{i,j-1}-d_{ij}|  &=&  O(i+1).
\end{eqnarray}
\eppsn

\bcrlre
Let $(\cla_f, \clh, D)$ be a spectral triple equivariant under the action of 
$S^1\times S^1$. Then $D$ can not be $p$-summable if $p<2$.
\ecrlre
\prf
This is a consequence of  the following growth
restriction on the $d_{ij}$'s:
\begin{equation}
d_{ij}=O(i+|j|+1),\label{eq:willsee}
\end{equation}
which follows from the last proposition.\qed

That there indeed exists a spectral triple that is 2-summable is easy to see,
by just taking $D$ to be the operator 
\be\label{genericd}
D=N\otimes S +I\otimes N,
\ee
where 
$S=\sum_{i\geq 0\atop j\geq 0} |e_{ij}\rangle\langle e_{ij}|-
    \sum_{i\geq 0\atop j< 0} |e_{ij}\rangle\langle e_{ij}|$.

\brmrk\rm
The obstruction element given by Voiculescu~(\cite{voi})
turns out to be zero for the ideals $\cll^{(p,\infty)}$ where $p<2$
(note that by proposition~1.7, \cite{voi}, it is enough to look at
positive finite-rank contractions from the commutant $U(G)^\prime$
in order to calculate this obstruction). Though one can not conclude anything definite from
this, it is possible that by dropping the condition of $S^1\times S^1$-equivariance,
Dirac operators of lower summability might be achievable.
\ermrk

\bppsn
Let $D$ be as in the previous proposition. Assume that $D$ has compact resolvent.
Then up to a compact perturbation, we have
\be\label{what}
\left.
\begin{minipage}{350pt}
\begin{enumerate}
\item For each $j\in\bbz$, all the $d_{ij}$'s are of the same sign,
 \item there is a big enough integer $M$ such that 
   \begin{enumerate}
   \item all the $d_{ij}$'s for $j\geq M$ are of the same sign, 
   \item all the $d_{ij}$'s for $j\leq -M$ are of the same sign.
     \end{enumerate}
\end{enumerate}
\end{minipage}
\right\}
\ee
\eppsn
\prf Again, the proof is very similar to the proof of
proposition~3.2 in \cite{CP}, and hence is omitted.\qed

This proposition says in particular that if $(\cla,\clh, D)$ is a $G$-equivariant
Fredholm module, then upto a compact perturbation, $P=\frac{I+\mbox{sign\,}D}{2}$ must be 
one of the following, where $E$ is some finite subset of 
$\{-M+1,-M+2,\ldots, M-1\}$:
\begin{eqnarray*}
P_1=\sum_{i\geq 0\atop j\leq -M} |e_{ij}\rangle\langle e_{ij}|
 +\sum_{i\geq0\atop j\in E} |e_{ij}\rangle\langle e_{ij}|,  &&
P_2=\sum_{i\geq 0\atop j\geq M} |e_{ij}\rangle\langle e_{ij}|
 +\sum_{i\geq0\atop j\in E} |e_{ij}\rangle\langle e_{ij}|,\\
P_3=\sum_{i\geq0\atop j\in E} |e_{ij}\rangle\langle e_{ij}|,&&
P_4=\sum_{i\geq0\atop j\in E^c} |e_{ij}\rangle\langle e_{ij}|.
\end{eqnarray*}

We will prove below that the $D$ given by (\ref{genericd}) is
in some sense the unique nontrivial Dirac operator for the representation
$\pi$ of $\cla$.

\bthm
Let $D'$ be an $G$-equivariant Dirac operator. Then the Kasparov module 
associated with $D'$ is either trivial or is same as the one 
associated with $D$ or $-D$.
\ethm
 \prf
Let $u=\chi_{\{0\}}(\beta^*\beta)(\beta-I)+I$.
First, observe that $\langle [u],(\cla,\clh,D)\rangle=\mbox{index}\, SuS=1$.
Since the $K$-groups for $SU_q(2)$ are free abelain, by the results of
Rosenberg \& Schochet (\cite{RS}), it is now enough to
show that $\langle [u],(\cla,\clh,D')\rangle$ is either 0 or $\pm 1$
if $P':=\frac{I+\mbox{sign\,}D'}{2}$ is one of the $P_i$'s above.
Since $\langle [u],(\cla,\clh,D')\rangle=\mbox{index}\, P'uP'$,
direct calculation now  tells us that if $P'$
is $P_3$ or $P_4$, the above pairing would be zero; it would be $-1$ if $P'=P_1$, 
and it is $1$ and if $P'=P_2$.\qed\\
The canonical unitary $\left(\brray{lr}\alpha &-q\beta^*\cr \beta &\alpha^*\erray\right)$
that comes in the definition of $SU_q(2)$ has non-trivial $K$-theory class
(see the remark following theorem~5, \cite{C8}). One can verify that
by computing its pairing with $D\otimes I$ (acting on $\clh\otimes\bbc^2$).

The following proposition can be derived as a corollary to 
proposition~4.3, \cite{CP}. But the proof presented there 
was just by computing pairings between appropriate
elements and does not give an insight as to why it is true.
We give a different proof here that sheds light on this.

\bppsn
Given any $m\in K^1(SU_q(2))=\bbz$, there exists a Kasparov module
$(L_2(h), F)$ which induces this element.
\eppsn
\prf
Using remark~\ref{interpret}, one could look at $L_2(h)$ as
$\oplus \cla p_i$. Representation of $\cla$ by left multiplications
in each piece looks like $\pi$. Now given $m$ in $\bbz$,
one has to pick $m$ copies of $\pi$, and define $F$ to be 
$(\mbox{sign}\,m)S$ on each of these pieces and $I$ on others. Then 
$(L_2(h), F)$ would be the required module.\qed


\newsection{Connes-de Rham complex}\label{complex}
Let $\Omega^\bullet(\cla_f)=\oplus_n\Omega^n(\cla_f)$ be the 
universal graded differential algebra over $\cla_f$,
i.e.\
$\Omega^n(\cla_f)=\mbox{span}\{a_0(\delta a_1)\ldots 
 (\delta a_n): a_i\in\cla_f, \delta(ab)=a(\delta b)+(\delta a)b\}$.
The universal differential algebra is not very interesting from the
cohomological point of view. Interesting cohomologies are obtained
from the representations of the algebra.
For the spectral triple $(\cla_f,\clh,D)$, one has the standard 
Connes-de Rham complex of noncommutative exterior forms 
$\Omega^\bullet_D(\cla_f)$, given by
\[
\Omega^\bullet_D(\cla) :=
\Omega^\bullet(\cla)/(\mathfrak{K}+\delta \mathfrak{K}) \cong
\pi(\Omega^\bullet(\cla))/\pi(\delta \mathfrak{K}).
\]
where $\mathfrak{K} = \oplus_{p \ge 0}\mathfrak{K}_p$ is the two sided ideal 
of $\Omega^\bullet(\cla)$ given by 
$\mathfrak{K}_p=\{\omega \in \Omega^p(\cla) : \pi (\omega) =0 \}$.
But often, the explicit computation of this complex
is rather difficult. What we will do is the following.
We will compute the complex
obtained from the representation $\theta \circ \pi: \Omega^\bullet
(\cla) \raro \clq (\clh)$ where 
$\theta : \cll(\clh) \raro \clq (\clh)= \cll(\clh)/\clk(\clh)$ 
is the projection onto the Calkin algebra.  
More specifically, let $\widetilde{d}:\cla_f\raro \cll(\clh)$   
be given by
$\widetilde{d}a=[D,\pi(a)]$.
Define $\pi_n:\Omega^n(\cla_f)\raro \cll(\clh)$ by
$\pi_n(a_0(\delta a_1)\ldots (\delta a_n))=\pi(a_0)(\widetilde{d}a_1)\ldots(\widetilde{d}a_n)$.
Define $d=\theta\circ\widetilde{d}$,
$\psi_n=\theta\circ\pi_n$, and $\psi:=\oplus\psi_n:\oplus\Omega^n\raro \clq(\clh)$.
Let $J_n=\ker \psi_n$.
Define
$\Omega_d^n(\cla_f)=\Omega^n(\cla_f)/(J_n+\delta J_{n-1})$.

Then $\Omega_d^n(\cla_f)=\psi(\Omega^n(\cla_f))/\psi(\delta J_n)$.
We will compute these cohomologies $\Omega_d^n(\cla_f)$. Before  entering
the computations, it should be stressed here that by computing these 
rather than the standard complex, we do not lose much. Because, first, 
since for a compact operator $K$ one has
$\mbox{Tr}_\omega (K|D|^{-2})=0$, proposition~5, page~550,  \cite{C7}
concerning the Yang-Mills functional holds in our present case.
Second, in the context of the canonical spectral triple 
associated with a compact Riemannian spin manifold this prescription 
also gives back the exterior complex.

First, we need the following lemma which will be very useful for the computations.

\blmma \label{tech1}
Assume $a,b\in\cla_f$ and $c\in\clk(\clh)$. If $a(I\otimes S)+b=c$, then
$a=b=0$.
\elmma
\prf
For a functional $\rho$ on $\cll(L_2(\bbn))$, and $T\in\cll(\clh)$,
denote by $a_\rho$ the operator $(\rho\otimes \mbox{id})T$.
Now observe that for any $a\in\cla_f$ and any functional $\rho$,
\be \label{eq:comm2}
a_\rho\ell =\ell a_\rho.
\ee

Write $P=\frac{1}{2}(I+S)$. It is easy to see that the given condition
implies that
$(b_\rho-a_\rho) + 2a_\rho P=c_\rho$,
which in turn implies that
\bea
(b_\rho-a_\rho)e_i  &=& c_\rho e_i \quad\forall i<0,\label{eq:temp1}\\
(b_\rho+a_\rho)e_i &=& c_\rho e_i \quad\forall i\geq0.\label{eq:temp2}
\eea

Now from (\ref{eq:comm2}) and (\ref{eq:temp1}), it follows that for any
$i,j\in\IZ$ and $j<0$,
\bean
\|(b_\rho-a_\rho)e_i\| &=& \|(b_\rho-a_\rho)\ell^{j-i}e_j\|\\
&=&\|\ell^{j-i}(b_\rho-a_\rho)e_j\| \\
&=& \|(b_\rho-a_\rho)e_j\|  \\
&=& \|c_\rho e_j\|.
\eean
Since $c$ is compact, $\lim_{j\raro -\infty}\|c_\rho e_j\|=0$.
Hence $(b_\rho-a_\rho)e_i=0$ for all $i$. In other words,  $(b_\rho-a_\rho)=0$.
Since this is true for any $\rho$,
we get $a=b$.
Using this equality, together with equations (\ref{eq:comm2}) and (\ref{eq:temp2}), a similar
reasoning yields $a=0$.
\qed

\blmma \label{cohom1}
Let $\cli_\beta$ denote the ideal in $\cla_f$ generated by $\beta$ and $\beta^*$.
Then for $n\geq 1$, we have
\begin{equation} \label{eq:cohom1}
\psi(\Omega^n(\cla_f))=(I\otimes S)^n\cla_f + (I\otimes S)^{n+1}\cli_\beta.
\end{equation}
\elmma
\prf
Let us first prove the equality for $n=1$.
Let $Z_k=q^{N+k}(N+k)$,
$B_{jk}=\sum_{i=j-k+1}^j |e_{i-1}\rangle\langle e_i|$, and
\[
C_j=\cases{\sum_{i=0}^{j-1}|e_i\rangle\langle e_{i-1}| & if $j\geq 1$,\cr
                                    0 &if $j=0$.}
\]
It is easy to check that
\bea
[D,\alpha]  & = & \alpha (-I\ots S),\nn \\
{}[D,\beta] & =&  q^N N\ots [S,\ell^*] +\beta.
\eea
It follows from these that
\bea \label{eq:comm3}
[D, \alpha_i\beta^j{\beta^*}^k] &=& -i(I\otimes S)\alpha_i\beta^j{\beta^*}^k
      +(j-k)\alpha_i\beta^j{\beta^*}^k
         + 2(Z_i\otimes C_j)\alpha_i\beta^{j-1}{\beta^*}^k \nn \\
       && \qquad \qquad     -2(Z_i\otimes B_{jk})\alpha_i\beta^j{\beta^*}^{k-1}.
\eea
Hence $d(\alpha_i\beta^j{\beta^*}^k)=-i(I\otimes S)\alpha_i\beta^j{\beta^*}^k
      +(j-k)\alpha_i\beta^j{\beta^*}^k$.
Thus for any $a\in\cla_f$,
\be \label{eq:temp3}
da = (I\otimes S)b+c, \qquad\mbox{ where } b\in\cla_f, \quad c\in\cli_\beta.
\ee
Note that for any $a'\in\cla_f$, $\psi(a')(I\otimes S)=(I\otimes S)\psi(a')$ in $\clq(\clh)$.
Hence
$\psi(a'(\delta a))$ is again of the form $(I\otimes S)b+c$, where $b\in\cla_f$, $c\in\cli_\beta$,
i.e.\ is a member of $(I\otimes S)\cla_f + \cli_\beta$.
Thus $\psi(\Omega^1(\cla_f))\seq (I\otimes S)\cla_f + \cli_\beta$.
For the reverse inclusion, observe that
$(I\otimes S)= (1-q^2)^{-1}((d\alpha)\alpha^*+q^2 (d\alpha^*)\alpha)$,
$\beta=d\beta$ and $\beta^*=-d\beta^*$.

The inductive step follows easily from (\ref{eq:temp3}).
\qed
\blmma \label{cohom2}
$J_0=\{0\}$, and for $n\geq 1$, we have
\begin{equation}\label{eq:cohom2}
\psi(\delta J_n)=(I\otimes S)^{n+1}\cla_f + (I\otimes S)^{n+2}\cli_\beta.
\end{equation}
\elmma
\prf
By lemma~\ref{tech1}, $\psi:\cla_f\raro\clq(\clh)$ is faithful.
Hence it follows that $J_0=\{0\}$.

We will prove here (\ref{eq:cohom2}) by induction.
From lemma~\ref{cohom1}, we have
$\psi(\delta J_1)\seq \psi(\Omega^2(\cla_f))=\cla_f+(I\otimes S)\cli_\beta$.
Let us show that  $I$, $(I\otimes S)\beta$ and $(I\otimes S)\beta^*$
are all members of $\psi(\delta J_1)$.

Choose $\omega\in\Omega^1(\cla_f)$ such that $\psi(\omega)=(I\otimes S)$.
Let $\omega_k=k\alpha_k\omega - \delta(\alpha_k), k=\pm 1$.
Then it follows from (\ref{eq:comm1}) that
$\psi(\omega_k)=k\alpha_k(I\otimes S) - k\alpha_k(I\otimes S) =0$,
so that $\omega_k\in J_1$.
$\psi(\delta\omega_k)=\psi(k(\delta\alpha_k)\omega)=k^2\alpha_k=\alpha_k\in \psi(\delta J_1)$,
i.e.\  both $\alpha$ and $\alpha^*$ are in $\psi(\delta J_1)$.
It follows from this that $I\in \psi(\delta J_1)$.

Next we show that $(I\otimes S)\beta\in\psi(\delta J_1)$.
Take $\omega=\frac{1}{2}(\alpha(\delta\beta) -\delta(\alpha\beta) +q\beta(\delta\alpha))$.
Then $\psi(\omega)=0$ and $\psi(\delta\omega)=(I\otimes S)\alpha\beta$.
So
$(I\otimes S)\alpha\beta\in \psi(\delta J_1)$.
Similarly taking
$\omega=\frac{1}{2}(\alpha^*(\delta\beta) -\delta(\alpha^*\beta) +q^{-1}\beta(\delta\alpha^*))$,
it follows that
$(I\otimes S)\alpha^*\beta\in \psi(\delta J_1)$.
These two together imply $(I\otimes S)\beta\in \psi(\delta J_1)$.

A similar argument shows that $(I\otimes S)\beta^*$ is also in $\psi(\delta J_1)$.
Thus $\cla_f+(I\otimes S)\cli_\beta=\psi(\delta J_1$).

For the inductive step, notice that
$\psi(\delta J_n)\seq \psi(\Omega^{n+1}(\cla_f))=(I\otimes S)^{n+1}\cla_f + (I\otimes S)^{n+2}\cli_\beta$.
We will show that the following are all
elements of $\psi(\delta J_n)$:
\[
\begin{array}{lll}
(I\otimes S)^{n+1}\alpha,  &  (I\otimes S)^{n+2}\alpha\beta,  &  (I\otimes S)^{n+2}\alpha\beta^*,\\
(I\otimes S)^{n+1}\alpha^*,  &  (I\otimes S)^{n+2}\alpha^*\beta,  &   (I\otimes S)^{n+2}\alpha^*\beta^*.
\end{array}
\]
From the right $\cla_f$-module structure of $\psi(\delta J_n)$, it will then follow that
$(I\otimes S)^{n+1}$, $(I\otimes S)^{n+2}\beta$ and $(I\otimes S)^{n+2}\beta^*$ are in
$\psi(\delta J_n)$, giving us the other inclusion.

Choose $\omega\in J_{n-1}$ such that $\psi(\delta\omega)=(I\otimes S)^n$.
Take $\omega_k=k\omega(\delta\alpha_k)$, $k=\pm 1$.
Then $\omega_k\in J_n$ and
$\psi(\delta\omega_k)=(I\otimes S)^{n+1}\alpha_k$.
Similarly choosing $\omega$ such that $\psi(\delta\omega)=(I\otimes S)^{n+1}\beta$
and $\omega_k$ as before, we get
$\omega_k\in J_n$ and
$\psi(\delta\omega_k)=q^{-k}(I\otimes S)^{n+2}\alpha\beta$.
Finally, take $\omega$ such that $\psi(\delta\omega)=(I\otimes S)^{n+1}\beta^*$
and $\omega_k$ as before to show that
$(I\otimes S)^{n+2}\alpha_k\beta^*\in\psi(\delta J_n)$.
\qed
\bthm \label{cohom3}
\begin{displaymath}
\Omega_d^n(\cla_f)=\cases{\cla_f\oplus\cli_\beta & if $n=1$,\cr
                                             \{0\} & if $n\geq 2$.}
\end{displaymath}
\ethm
\prf
Proof follows from lemmas \ref{cohom1} and \ref{cohom2}.
\qed

\brmrk \rm
The differential $d:\cla_f\rightarrow\Omega_d^1(\cla_f)=\cla_f\oplus\cli_\beta$
is given by
\[
d(\alpha_i\beta^j\beta^{*k})
 =-i\alpha_i\beta^j\beta^{*k}\oplus (j-k)\alpha_i\beta^j\beta^{*k}.
\]
\ermrk
\section{$L^2$-complex of Frohlich et.\ al.}
In this section we will compute the complex of square 
integrable forms for the spectral triple $(\cla_f,\clh,D)$.
For that we begin with similar computations for the spectral 
triple  $(\IC [ z ,z^{-1}],\clh_0=L_2(\IZ),D_0=N)$ associated 
with the algebra $\IC[z,z^{-1}]$.
Here we consider the embedding $\pi_0:  \IC[z,z^{-1}] \raro \cll(\clh)$ 
that maps $z$ to $\ell$.
\blmma
{\rm (i)} $\widetilde{\Omega}^n_{D_0}(\IC[z,z^{-1}])  =  0 , \; \mbox{for }n \ge 2, $\\
{\rm (ii)} $\widetilde{\Omega}^1_{D_0}(\IC[z,z^{-1}])  =  \IC[z,z^{-1}].$
\elmma
\prf (i) Let $\omega =\sum \alpha_{n_0,\cdots,n_k}z^{n_0} \delta z^{n_1} \cdots 
\delta z^{n_k} \in \Omega^k(\IC[z,z^{-1}])$, where the sum is a finite one
and $\delta$ is the universal differential. Then it is easily verified that 
\[
( \omega,\omega)_{D_0} = \int {( \sum n_1 \cdots n_k \alpha_{n_0,\cdots,n_k}  z^{\sum_0^k n_j})}^* 
( \sum n_1 \cdots n_k \alpha_{n_0,\cdots,n_k}  z^{\sum_0^k n_j}) dz,
\]
where $dz$ is the Lebesgue measure on the circle.
Therefore,
\bean
\mathfrak{K}_k(\IC[z,z^{-1}])& :=& \{ \omega \in \Omega^k(\IC[z,z^{-1}]) : \; ( \omega,\omega)_{D_0} =0 \} \cr
& =& \{ \sum \alpha_{n_0,\cdots,n_k}z^{n_0} \delta z^{n_1} \cdots 
\delta z^{n_k}  : \; \sum_{n_0+\cdots +n_k=r} n_1 \cdots n_k \alpha_{n_0,\cdots,n_k} =0, \forall r \} . 
\eean
Consequently we have,
\bea
z^{n_0} \delta z^{n_1} \cdots 
\delta z^{n_k}-n_1 \cdots n_k z^{\sum_0^kn_i-k} 
  \delta z \cdots \delta z & \in & \mathfrak{K}_k(\IC[z,z^{-1}]), \label{S:lem:5.15}\\
\delta z^{r}   \delta z\cdots \delta z - r z^r \delta z\cdots \delta z & 
   \in & \mathfrak{K}_k(\IC[z,z^{-1}]), \label{S:lem:5.16}\\
z^r \delta z\cdots \delta z -\frac{1}{r+1} \delta z^{r+1}   
\delta z\cdots \delta z & \in  & \mathfrak{K}_{k-1}(\IC[z,z^{-1}]). \label{S:lem:5.17}
\eea
From~(\ref{S:lem:5.17}) we get $\delta z^{r}   \delta z\cdots \delta z  \in   \delta \mathfrak{K}_{k-1}(\IC[z,z^{-1}])$. Combining this with (\ref{S:lem:5.15}) and (\ref{S:lem:5.16}) we get, 
\[
z^{n_0} \delta z^{n_1} \cdots 
\delta z^{n_k} \in \mathfrak{K}_k(\IC[z,z^{-1}])
  +\delta \mathfrak{K}_{k-1}(\IC[z,z^{-1}]) \mbox{ for large } n_0. 
\]
Since $ \mathfrak{K}_k(\IC[z,z^{-1}])+\delta \mathfrak{K}_{k-1}(\IC[z,z^{-1}])$ is a  bimodule we have
\[
z^{n_0} \delta z^{n_1} \cdots 
\delta z^{n_k} \in \mathfrak{K}_k(\IC[z,z^{-1}])+\delta \mathfrak{K}_{k-1}(\IC[z,z^{-1}]) 
  \quad \forall\; n_0,\cdots,n_k.
\]
This proves (i).\\
(ii) It suffices to note that
\[
z^{n_0} \delta z^{n_1} - n_1 z^{n_0+n_1-1} \delta z \in \mathfrak{K}_1(\IC[z,z^{-1}]).
\]
The induced 
$ d: {\widetilde{\Omega}}^0_{D_0} (\IC[z,z^{-1}]) \raro \IC[z,z^{-1}]$ 
is given by $d(z^n)= n z^n $. \qed 

Now we are in a position to compute the complex of square integrable 
forms.
\bthm
{\rm (i)} $ {\widetilde {\Omega }}^n_D ( \cla_f) =0 $ for $n \ge 2 $. \\
{\rm (ii)} $ {\widetilde {\Omega }}^n_D ( \cla_f) = \IC[z,z^{-1}]$ for $n=0,1$ 
(equality as an $\cla_f$ bimodule), and the differential 
$d:\cla_f\rightarrow {\widetilde {\Omega }}^1_D ( \cla_f)$
is given by $d(\alpha_i\beta^j\beta^{*k})=-iz^i$.
\ethm
\prf Note that the homomorphism $\sigma$ in (\ref{100}) 
induces a surjective homomorphism denoted by the same symbol from 
$\cla_f$ to $\IC[z,z^{-1}]$. We have the following
 short exact sequence
\[
0\longrightarrow \cli_\beta\longrightarrow\cla_f
\stackrel\sigma\longrightarrow\IC[z,z^{-1}]\longrightarrow 0,
\]
Let $\sigma_k: \Omega^k(\cla_f) \raro \Omega (\IC[z,z^{-1}])$ be the 
induced surjective map. One easily verifies that 
$(\omega,\omega)_D=(\sigma_k(\omega),\sigma_k(\omega))_{D_0}$. Therefore,
\[
\mathfrak{K}_k(\cla_f)=\{ \omega \in \Omega^k (\cla_f) : \; 
(\omega, \omega)_D=0 \} = \sigma_k^{-1}(\mathfrak{K}_k( \IC [ z,z^{-1}])). 
\]
We have the following commutative diagram
\newcommand{\bbcz}{\mathbb{C}[z,z^{-1}]}

\begin{center}
\begin{tabular}{ccccccc}
$\mathfrak{K}_0=I_\beta$& $\longrightarrow$ & $\mathcal{A}_f$ &$\stackrel{\sigma}\longrightarrow$& $\bbcz$ &
        $\longrightarrow$ & $\pi_0(\bbcz)$\\
&&$\downarrow$&&$\downarrow$&&  $\downarrow$\\
$\mathfrak{K}_1(\cla_f)$& $\longrightarrow$ & $\Omega^1(\mathcal{A}_f)$ &$\stackrel{\sigma_1}\longrightarrow$& $\Omega^1(\bbcz)$ &
        $\longrightarrow$ & $\tilde{\Omega}_{D_0}^1(\bbcz)$\\
&&$\downarrow$&&$\downarrow$&&  $\downarrow$\\
$\mathfrak{K}_2(\cla_f)$& $\longrightarrow$ & $\Omega^2(\mathcal{A}_f)$ &$\stackrel{\sigma_2}\longrightarrow$& $\Omega^2(\bbcz)$ &
        $\longrightarrow$ & $\tilde{\Omega}_{D_0}^2(\bbcz)$\\
&$\ldots$&&$\ldots$&&$\ldots$& $\ldots$\\
$\mathfrak{K}_n(\cla_f)$& $\longrightarrow$ & $\Omega^n(\mathcal{A}_f)$ &$\stackrel{\sigma_n}\longrightarrow$& $\Omega^n(\bbcz)$ &
        $\longrightarrow$ & $\tilde{\Omega}_{D_0}^n(\bbcz).$
\end{tabular}
\end{center}
This along with the previous lemma proves the theorem. 
We will only illustrate (i). 

Let $\omega_n \in \Omega^n(\cla_f)$, then by the previous lemma 
$\sigma_n(\omega_n)= \omega_{1,n}+\delta\omega_{2,n-1}$ where 
$\omega_{1,n}\in \mathfrak{K}_n(\IC[z,z^{-1}]),\omega_{2,n-1}\in \mathfrak{K}_{n-1}(\IC[z,z^{-1}])$. 
Let $\omega_{1,n}^\prime=\sigma_n^{-1}(\omega_{1,n}),
\omega_{2,n-1}^\prime=\sigma_{n-1}^{-1}(\omega_{2,n-1})$, 
then $\sigma_n(\omega_n-\omega_{1,n}^\prime-\delta \omega_{2,n-1}^\prime)=0$ 
implying $\omega_n \in \mathfrak{K}_n+\delta \mathfrak{K}_{n-1}$. \qed

\section {Computations for the quantum sphere}
In this section we will briefly indicate how to carry out 
all the earlier constructions for the quantum spheres.
We will be sketchy because  most of the arguments 
are very similar to the case of $SU_q(2)$.
Quantum sphere was introduced by Podle\'s in \cite{PO}. 
This is the universal C*-algebra, denoted by $C(S^2_{qc})$, generated 
by two elements $A$ and $B$ subject to the following relations:
\bean
A^* = A ,&&  B^*B = A - A^2 + c I ,\\
BA = q^2 AB,&&BB^* = q^2 A - q ^ 4 + c I.
\eean
Here the deformation parameters $ q$ and $c$ satisfy $ | q | < 1, c>0 $. 
For later purpose we also note down two irreducible representations 
whose direct sum is faithful. Let 
$\clh_{+}  = l^2 ( \bbn  ), \clh_{-} = \clh_{+}$. 
Define $ \pi_{\pm} (A), \pi_{\pm} (B): \clh_\pm \raro \clh_\pm$ by 
\bean
\pi_{\pm}(A) ( e_n)= \lambda_\pm q^{2n} e_n \;\;\quad \; &
 {\rm   where   }& \quad \; \lambda_\pm = \frac {1}{ 2}  
       \pm {(c + \frac {1}{4} )}^{1/2}, \cr
\pi_{\pm}(B) ( e_n)= {c_\pm (n)}^{1/2} e_{n-1} &{\rm  where}&
 c_\pm (n)= \lambda_\pm q^{2n}- {(\lambda_\pm q^{2n})}^2 + c.
\eean
 Since $\pi =\pi_+\oplus \pi_-$ is a faithful representation,  an 
immediate corollary follows.
\bthm [Sheu \cite{S}]
\emph{(i)} $C (S^2_{qc}) \cong C^* (\mathscr{T})\oplus_\sigma 
  C^*(\mathscr{T}) := \{(x,y) : \; x,y \in 
  C^*(\mathscr{T}), \sigma(x)=\sigma(y) \} $ 
where $C^*(\mathscr{T})$ is the Toeplitz algebra and 
$\sigma : C^*(\mathscr{T}) \raro C(S^1)$ is the symbol homomorphism.\\
\emph{(ii)} We have a short exact sequence 
\bea
\label{S:ses:2}
0\longrightarrow \clk \stackrel i\longrightarrow  C(S^2_{qc}) 
\stackrel\alpha \longrightarrow C^*(\mathscr {T}) \longrightarrow 0
\eea
\ethm
\prf (i) An explicit isomorphism is given by $x \mapsto (\pi_+(x),\pi_-(x))$. \\
(ii) Define $\alpha((x,y))=x$ then $\ker \alpha =\clk$. \qed
\bcrlre
\emph{(i)} $ K_0(C(S^2_{qc}))=K^0(C(S^2_{qc}))=\IZ \oplus \IZ.$  \\
\emph{(ii)} $ K_1(C(S^2_{qc}))=K^1(C(S^2_{qc}))=0.$
\ecrlre
\prf The six term exact sequence associated with (\ref{S:ses:2}) 
along with the KK-equivalence of $\clk$ and $C^*(\mathscr{T})$ 
with $\IC$ proves the result.\qed

\bppsn
Let $\cla_{fin}$ be the *-subalgebra of $C(S^2_{qc})$ generated by $A$ and $B$.
Then 
\[
\left(\cla_{fin},\;\clh=\clh_+\oplus \clh_-,\; D= {0\;\;N \choose N\;\; 0},
\;\gamma=\left({ 1\atop 0}{0 \atop {-1}}\right)\right)
\] 
is an even spectral triple.
\eppsn 
\prf We only have to show that $[D,a]$ is bounded for 
$a \in \cla_{fin}$. For that it is enough to note that\\
(i) $N \pi_\pm(A),\; \pi_\pm (A) N$ are bounded,\\
(ii) $n ({c_\pm(n)}^{1/2} - \sqrt{c} )$ is bounded as $n$ becomes large,\\
(iii) $[N,l]=l$. \qed
\brmrk  \rm
This spectral triple has nontrivial Chern character. 
This can be seen as follows: let 
$P_0=i(|e_0 \rangle \langle e_0|) \in C(S^2_{qc})$, 
then applying proposition~4, page 296, \cite{C7},
we get the 
index pairing $\langle [P_0],[(\cla_{fin},\clh,D,\gamma)]\rangle=-1$, 
implying nontriviality of the spectral triple.
\ermrk
Now we will briefly indicate the computations of the 
complex $(\Omega^\bullet_d(\cla_{fin}),d)$ introduced 
at the beginning of section~\ref{complex}.
\bthm
\emph{(i)} $\Omega^n_d(\cla_{fin})=0$ for $n \ge 2 $.\\
\emph{(ii)} $\Omega^1_d(\cla_{fin})= \IC[z , z^{-1}]$, 
here also equality is as an $\cla_{fin}$ bimodule.
\ethm
\prf Let $\pi$ be the associated representation of 
$\Omega^\bullet(\cla_{fin})$ in $\cll(\clh)$. 
Then straightforward verification gives (i) $[D,A]$ is compact, 
(ii) $[D,B]= l \otimes \kappa + \mbox{compact}$, and 
(iii) $[D,B^*]= -l^* \otimes \kappa +\mbox{compact}$, 
where $\kappa=\left( \matrix{ 0 & 1 \cr 1 & 0 \cr }\right)$. 
Therefore,
modulo compacts 
\bean
\pi (\Omega^{2k+1}(\cla_{fin})) &=& C^*_{fin}(\mathscr{T})\otimes \kappa \cr
\pi (\Omega^{2k}(\cla_{fin})) &=& C^*_{fin}(\mathscr{T})\otimes I_2,
\eean
where $C^*_{fin}(\mathscr{T})$ is the *-algebra generated by $\mathscr T$.
Now for~(i), note that 
\[
\omega_n=B \delta B^* \underbrace{\delta B \cdots \delta B}_{n-2 \; times } 
+ B^* \delta B \underbrace{\delta B \cdots \delta B}_{n-2 \; times }
\] 
satisfies (a) $\pi(\omega_n)$ is compact and (b)$\pi(\delta \omega_n)=2I$ 
is invertible, hence (i) follows. \\
For~(ii),  observe that if $a \in \cla_{fin}$ and $\pi(a)$ is compact 
then $Na$ and $aN$ both compact. Hence, 
$\Omega^1_d(\cla_{fin})=\pi (\Omega^{1}(\cla_{fin}))=\IC[z,z^{-1}]$ 
because modulo compacts $\IC[z,z^{-1}]$ is $C^*(\mathscr{T})$. \qed


\noindent{\sc Partha Sarathi Chakraborty}\\
{\footnotesize Indian Statistical 
Institute,  203, B. T. Road, Calcutta--700\,108, INDIA\\[-.5ex]
         email: parthasc\_r@isical.ac.in}\\[1ex]
{\sc Arupkumar Pal}\\
         {\footnotesize Indian Statistical 
Institute, 7, SJSS Marg, New Delhi--110\,016, INDIA\\[-.5ex]
         email: arup@isid.ac.in}

\end{document}